\documentclass[12pt]{article}
\usepackage{amsmath,amssymb,amsthm,url}
\def\Z{{\mathbb Z}}  \def\R{{\mathbb R}}

\def\Cl{\mathop{\fam0 Cl}}

\long\def\comment#1\endcomment{}

\begin{document}


\centerline{\bf A short exposition of the Levine-Lidman example of spineless 4-manifolds\footnote{This is an extended version of math review of [LL].
I would like to acknowledge S. Akbulut, R. Karasev, A. Levine, T. Lidman and A. Zhubr for useful discussions.}}

\smallskip
\centerline{\bf A. Skopenkov\footnote{\texttt{https://users.mccme.ru/skopenko}, Moscow Institute of Physics and Technology, Independent University of Moscow.
Supported in part by the Russian Foundation for Basic Research Grant No. 19-01-00169 and by Simons-IUM Fellowship.}}

\smallskip
The paper [LL] outlines a proof of an interesting result in topology of manifolds.
Here we present a shorter (and hopefully clearer) exposition, see Remark 4.a.
We reveal that some parts of the proof are missing, and that some results are used in [LL] without proof or reference, or even without explicit statement.
Since the authors refused to make necessary details publicly available (see Remark 5), we invite a reader to publish them, thus completing the proof.

\smallskip
{\bf Conjecture 1.} [LL] {\it There is a compact smooth 4-manifold $W$ with boundary such that $W$ is homotopy equivalent to $S^2$ but there does not exist an embedding $S^2\to W$ which is a homotopy equivalence and is
simplicial for some triangulations of $W$ and of $S^2$.}\footnote{We do not formally use the notion of {\it spine} because it is defined differently in different areas of topology of manifolds, and because [LL, Theorem 1] needs an explanation why spine is defined for a PL manifold but is used for smooth manifold.}

\smallskip
The well-known Kirby's list has this as a problem motivated by generalizations of the celebrated Browder-Casson-Haefliger-Sullivan-Wall Theorem of 1960s.
See a short and clear exposition of the background in [LL, the paragraph after Theorem 1.1]

The statement of [LL, Theorem 1.1] is more general.
The outline of the proof exposed below will presumably work for the generalization.

Perhaps it was non-trivial to invent the example of Conjecture 1.
The idea of proof is an easy application of Conjecture 2 below, which presumably is easily proved using recent known results on Heegard-Floer $d$-invariants (Conjecture 3 below).

Denote by $Q$ the total space of a circle bundle over $\R P^2$ with normal Euler number $-7$, and by $\Sigma$ the Seifert homology sphere $\Sigma(2,3,7)$.
A {\bf homology cobordism} between two oriented 3-manifolds $M_0,M_1$ is an oriented cobordism $W$ from $M_0$ to $M_1$ such that the inclusion of either $M_i$ induces an isomorphism on homology with integer coefficients.

\smallskip
{\bf Conjecture 2.} {\it The manifold $Q\#\Sigma$ is not homology cobordant to $\chi_L:=\chi_{L,4}(S^3)$ for any knot $L$.}

\smallskip
{\it Outline of proof of Conjecture 1 modulo Conjecture 2.}
Denote $p:=1$ in order to use figures from [LL].
Let $K$ be the meridian of the $p$-framed surgery curve as shown in framed surgery representation of $-\Sigma$
in [LL, Fig. 2].
(The knot $K$ is a singular fiber in a Seifert fibration of $-\Sigma$).
The 0-framing on $K$ (viewed as a knot in $S^3$) corresponds to the $+4$ framing on $K$ (as a knot in $-\Sigma$). Performing surgery using this framing produces $Q$, since we can cancel the $p$-framed component with its
0-framed meridian to produce [LL, Fig. 1(c)] with $m=−7$.
Let $W$ be the oriented 4-manifold obtained from $(-\Sigma\backslash B^3)\times[0,1]$ by attaching a $+4$-framed 2-handle along $K\times\{1\}$.
Since $\partial\left((-\Sigma\backslash B^3)\times[0,1]\right)\cong (-\Sigma)\#\Sigma$, we have
$\partial W\cong Q\# \Sigma$.

Clearly, homology groups of $(-\Sigma\backslash B^3)\times[0,1]$ are trivial.
Hence $W$ has the same homology as that of $S^2$.
Using the fact that 2,3,7 are pairwise relatively prime, [LL, Lemma 3.2] shows that $W$ is simply connected.
Hence $W$ is homotopy equivalent to $S^2$.

Assume to the contrary that there is an embedding $S^2\to W$ which is a homotopy equivalence and is
simplicial for some triangulations of $W$ and of $S^2$.
It is known and easy to see that the regular neighborhood $R$ of this embedding is obtained from the 4-ball by Dehn surgery along some knot $L$ with framing $+4$.
Clearly, $\Cl(W-R)$ is a cobordism between $\partial W$ and the boundary $\partial R=\chi_L$.
{\it We need to show that $\Cl(W-R)$ is a homology cobordism.}
Such facts are usually proved (or disproved) using Mayer-Vietoris sequence.
Then Conjecture 2 would give a contradiction.


\smallskip
{\bf Conjecture 3.} {\it To any closed oriented
3-manifold $M$ such that $H_1(M)$ is finite and any class
$s\in H_1(M)$ there corresponds an integer $\delta(M,s)$ such that the following holds.


(1) For each homology cobordism $Z$ between $M_0$ and $M_1$, and $s_j\in H_1(M_j)$ going to the same element of $H_1(Z)$ under the inclusions, $\delta(M_0,s_0)=\delta(M_1,s_1)$.

(2) Let $X$ be a smooth compact oriented 4-manifold (with boundary) having the same homology as that of $S^2$, and such that the self-intersection number of a generator $\alpha$ of $H_2(X)$ is $+4$.
Then $\delta(\partial X,\partial t)\equiv (t\cap\alpha)^2-1\mod8$ for any $t\in H_2(X,\partial)$.

(3) Let $X$ be the (properly oriented) trace of the surgery on any (oriented) knot in $S^3$ with framing $+4$, and $\alpha\in H_2(X)\cong\Z$, $\tau\in H_2(X,\partial)\cong\Z$ generators such that $\alpha\cap\tau=1$.
Then
$$|\delta(\partial X,\partial(t+\tau))-\delta(\partial X,\partial t)-2t\cap\alpha-1| \le8 \quad\text{for any}\quad t\in H_2(X,\partial).$$

(4) For
any $(M,s)$ we have $\delta(M\#\Sigma,s\#0) = \delta(M,s) + \delta(\Sigma,0)$.

(5) The number $\delta(\Sigma,0)$ is divisible by 8.

(6) The number $\delta(Q,s)$ assumes values -5,-9,0 and 0 for the four elements
$s\in H_1(Q)\cong\Z_4$.\footnote{In the statement of this property in [LL, p.6] one needs to delete `$\{$' and `$\}$' rather than to use them in an incorrect way (recall that $\{\dfrac{m+2}4,\dfrac{m-2}4,0,0\}=\{\dfrac{m+2}4,\dfrac{m-2}4,0\}$).}
}


\smallskip
Parts (2), (3), (6) are restatements of [OS03, Theorem 1.2], [HW16, Equation 2.3], [Doi15, \S3], respectively (see references in [LL]).
So these are results not conjectures.
(Part (3) is a restatement in a weaker form sufficient for the proof of Conjecture 2.)
{\it Part (5), (1), (4) were used in [LL] without reference, without explicit statement and even without mentioning, respectively.}

\comment

See reply from (October 22)

(1) follows from Theorem 1.2 in Ozsvath-Szabo's
"Absolutely graded Floer homologies and intersection forms for four-manifolds with boundary"

(4) follows from the same theorem, since being a group homomorphism makes the invariants additive.

(5)  The point is that the d-invariant of any integer homology sphere is an even integer.  There are several ways to see this, and it has appeared in various forms throughout the literature, although people tend not to give a proof because the ideas are standard in the area.  Here is a proof.
The d-invariant is defined to be the absolute grading of the bottom most element in the ``infinite tower'' in HF^+(Y,s).  The discussion above Proposition 3.3 in the same paper shows that the grading of this element has integral grading.  Now let's see that the grading is even.  Proposition 3.3 there shows that the mod 2 reduction of the grading lines up with the alternately-defined Z/2-grading in Heegaard Floer homology, coming from Lagrangian Floer homology constructions (i.e. the Z/2-grading coming from the sign of the intersections of the tori in the symmetric product).  That absolute grading has the property that the Euler characteristic of HF-hat is 1.  Using the long exact sequence between HF^+, HF^+, and HF-hat, the infinite tower contributes one dimension to HF-hat, and the torsion pieces (e.g things of the form F[U]/U^n) each contribute 2 dimensions to HF-hat in opposite parity of grading, one dimension is in the kernel of U and one is in the cokernel of U.  Therefore, the Euler characteristic is 1 if and only if the infinite tower is supported in Z/2-grading 0, i.e. even absolute grading.  Therefore, the d-invariant is an even integer.

A shorter, but less broad argument is as follows.  The Y_p are Seifert homology spheres, so their d-invariants can be computed using the formula in Corollary 1.5 of Ozsvath-Szabo's  "On the Floer homology of plumbed three-manifolds''.  Tere, K^2 + |G| is computing c_1^2(t) - 3\sigma(X) + 2\chi(X) where X is a negative-definite plumbing bounding the homology sphere Y_p (or possibly -Y_p, depending on the value of p).  Since c_1(t) is characteristic, c_1^2(t) - 3\sigma(X) + 2\chi(X) is a multiple of 8, and hence the d-invariant is an even integer.

\endcomment

\smallskip
{\it Outline of proof of Conjecture 2 modulo Conjecture 3.}
Suppose to the contrary that there is a homology cobordism $Z$ between $Q\#\Sigma$ and $\chi_L$.
Take $X,\alpha,\tau,t$ as defined in (3).
Then $\alpha\cap\alpha$ is $+4$ not $-4$.\footnote{In [LL] this was proved using Conjecture 3 in a more general situation of [LL, Proposition 3.3] unnecessary for the proof of Conjectures 1 or 2.}
Denote $X_+:=X\cup_{\partial X=\chi_L\subset\partial Z} Z$.
Denote by $\alpha_+$ the image of $\alpha$ under the inclusion-induced isomorphism $H_2(X)\to H_2(X_+)$.
Denote by $\tau_+,t_+$ the preimages of $\tau,t$ under the `cutting' isomorphism
$H_2(X_+,\partial)\to H_2(X,\partial)$.
{\it We need to show that $\partial t_+$ and $\partial t$, $\partial\tau_+$ and $\partial\tau$ go to the same element of $H_1(Z)$ under the inclusions.}
Such facts are usually proved (or disproved) using Mayer-Vietoris sequence.
Let us continue this outline assuming that this is proved.

By (2) the integers $\delta(Q\#\Sigma,i\partial\tau_+)$ are congruent modulo 8 to 7,0,3,0 for $i=0,1,2,3$, respectively.
By (4) and (5) we have $\delta(Q\#\Sigma,s\#0)\equiv \delta(Q,s)\mod8$ for any $s\in H_1(Q)$.
Hence by (4) and (6) the integers $\delta(Q\#\Sigma,k\partial\tau_+)-\delta(\Sigma,0)$ assume values -9,0,-5,0 for $k=0,1,2,3$, respectively.
Then
$$\delta(Q\#\Sigma,0)-\delta(Q\#\Sigma,3\partial\tau_+) = -9\not\in\{-1,7,15\}
\overset{(*)}\ni \delta(\chi_L,0)-\delta(\chi_L,3\partial\tau).$$
Here (*) holds by (3) applied to $t_+=3\tau_+$.
We obtain a contradiction to (1).

\smallskip
{\bf Remark 4.}
(a) Despite of being shorter than in [LL] the above exposition is not an alternative proof
but just a different presentation, making clear the structure and avoiding sophisticated language.
(I also do not repeat [LL, Fig. 1c, Fig. 2 and Lemma 3.2].)
When I use a specific theory, I explicitly state results (see Conjectures 2 and 3) to be proved with the help of this theory but in terms not involving the theory.
This makes the {\it application} of the result  accessible to mathematicians who have not specialized in the theory.
So one is motivated to study the theory and sees explicit statements which could guide this study.

Instead of the above way, some papers start exposition with details of a specific theory
which are matter-of-fact to specialists but are much less accessible to mathematicians from other (even close) areas.
This presumably happens because of the false assumption that the readers will accept artificial sophistication as depth and high non-triviality.


(b) The good old way of explicitly stating non-trivial results and either proving them or giving a reference is not only necessary for reliability.
Such a style contributes to the unity of mathematics, by making some areas more accessible to mathematicians from other areas.
This in turn ensures higher reliability.
Disregard of this style not only decreases reliability, but also contributes to artificial splitting of mathematics (and even of its areas) into different subjects whose representatives cannot use each other's work.

(c) This note contributes to countering a widespread opinion that mathematicians do not seriously check the work of each other, so that most of the published research is unreliable.

(d) In the notation of [LL] we have $\delta(M,s)=8d(M,\overline s)$, where the spin structure $\overline s$ extends to a spin structure $\overline t$ on a null-cobordism $X$ of $M$ as in Conjecture 3, (2) such that for any $t\in H_2(X,\partial X)$ with $\partial t=s$ we have $PDc_1(\overline t)\cap\alpha=t\cap\alpha$.
I multiplied the invariant by 8 so as not to bother a reader with (more complicated formulas involving) fractions.

Also note for comparison to [LL] that $n=4$, $i-2=t\cap\alpha$, $p=1$, $\Sigma=-Y_1$.

(e) Let us illustrate how much of Conjecture 3 can be recovered by standard algebraic topology.
Let $X$ be a smooth compact oriented 4-manifold (with boundary) having the same homology as that of $S^2$, and such that the self-intersection number of a generator $\alpha$ of $H_2(X)$ is 4.
Take any $t\in H_2(X,\partial)$.


$\bullet$ The number $t\cap\alpha$ depends on $t$ for given $(\partial X,\partial t)$.

Indeed, $\partial(t+4\tau)=\partial t$ and $(t+4\tau)\cap\alpha=t\cap\alpha+4$.

$\bullet$ The reduction modulo 4 of $t\cap\alpha$ only depends on $(\partial X,\partial t)$.

(Then $(t\cap\alpha)^2$ modulo 8 is independent of $t$ for given $(X,\partial t)$.)




Indeed, take two pairs $(X,t)$ and $(X',t')$ having the same boundaries $(M,s)$.
Let $U:=X\cup_Y X'$.
Since $H_1(M)\cong\Z_4$ is finite, we have $H_2(M)=0$.
Then by MVS
$$\ldots\to H_2(X)\oplus H_2(X')\to H_2(U)\to H_1(M)\to\ldots$$
the group $H_2(U)/\text{torsion}$ is generated by the images of $\alpha$ and $\alpha'$ under the inclusions $X\to U$ and $X'\to U$.
We denote these images by the same letters $\alpha$ and $\alpha'$.
Since $\partial t=\partial t'$, by another MVS
$$\ldots\to H_2(U)\overset{r\oplus r'}\to H_2(X,\partial)\oplus H_2(X',\partial)
\overset{\partial-\partial'}\to H_1(M)$$
there is a class in $H_2(U)$ going to $t\oplus t'$ under the homomorphism $r\oplus r'$ from the MVS.
This class is congruent modulo torsion to $p\alpha+p'\alpha'$ for some integers $p,p'$.
Then
$$t\cap\alpha-t'\cap\alpha' = (p\alpha+p'\alpha')\cap(\alpha-\alpha') = p\alpha^2-p'(\alpha')^2 = 4(p-p')$$
is divisible by 4.

$\bullet$ Denote by $\tau$ the generator of $H_1(X,\partial)\cong\Z$.
Then $((t+\tau)\cap\alpha)^2-(t\cap\alpha)^2=2\tau\cap\alpha+1$.

This follows because $\tau\cap\alpha=1$.


\smallskip
{\bf Remark 5.}  Here I present my letters to the authors of [LL].
These letters contain public statements of different reliability standards,
and my repeated invitations to the authors to add missing details described above.
These letters could be useful for other mathematicians (not necessarily math reviewers) who would like to politely invite
other authors to add missing details.\footnote{English is slightly corrected.
The letter of October 20 (and part of the letter of October 22) is omitted but are available upon request.
That letter describes yet another flaw to be corrected and repeats the invitation to update the arxiv version.}

\smallskip
{\it (October 14)} Dear Adam, Dear Tye,

It would be nice if you could update your interesting paper arxiv:1803.01765
adding the proof that `the $spin^c$ structures ... are identified through this cobordism'
in p. 4 at the end of proof of Theorem 2.2.
Then I will not need mention in my math review that this is not proved
(although simpler facts are proved).

Let me know if you are also interested to learn some remarks allowing to improve the exposition of your paper.
I would be glad not to mention them in my math review but just refer to arxiv update.

Best Regards, Arkadiy.

\smallskip
{\it (October 22)} Dear Adam, Dear Tye,

Thank you for your efforts on providing missing details for your arguments.
It would be a pity if these details will not be available to math community, so that
I would have to write in my review that the published proof is incomplete.
Hence it would be nice if you could update your interesting paper arxiv:1803.01765
adding explicit statements and either a proof or a reference for the results below that you use,
and  adding justification that (1') below can indeed be applied in p. 4 at the end of the proof of Theorem 2.2.

...

Best, A.

\smallskip
{\it (October 30)} Dear Adam, Dear Tye,

Attached please find a project of an extended version of Math Review to your interesting paper arxiv:1803.01765.
I would be most grateful for any remarks.
In particular, could you let me know if the additions marked with `???' are correct.

The project is written taking into account that you refused to update the arxiv version of your paper adding the required details as described above.
However, I would be glad to rewrite this note if the details will appear in the arxiv update of your paper
(together with a reference to an update of this note).

I do not mean that Conjectures 1, 2 and 3 are false, or that the details described above are hard to write.
I only mean that a published paper (and a Math Review) is for users, not for developers.
Working on details could be an interesting task for a developer but is usually not within intents of a user.
One of the best estimations of how hard details are is the amount of time required for authors
(or for other mathematicians) to make the details publicly available upon request of a reviewer.

Best Regards, Arkadiy.

\smallskip
{\it (October 30)} Dear Adam, Dear Tye,

Here's my reaction to Tye's letter which I just noticed and which perhaps was not intended to my eyes.

I am ready to publicly state that using a non-trivial result without explicit statement (and so without proof or reference)
makes the proof incomplete (in the sense of arxiv:1702.04259v1, text after Conjecture 2).
I would be glad if you publicly state the opposite.

But maybe a compromise would be more pleasant and effective.
The corrections required in your paper presumably are minor and very far from requiring an erratum
(I can only finally confirm this and the following judgments when I see the updated version.)
So your wish that the update should not be considered as an indication that the published version contains any serious gaps is presumably very much justified.
If you could update the arxiv version, I am willing to make this clear in my math review.
E.g. I can ignore the (presumably minor) gaps and just refer to the published paper together with the arxiv version.
Or I can explicitly write that `the arxiv version contains some minor details omitted in the published version'.
Alternatively, I would be glad to use a phrase you suggest.
I will also be glad to change (or to suppress) some less important criticism in the project of a note I sent you.

In fact, updating the arxiv version is not considered as an indication that the published version has any serious gaps.
Some authors previously updated arxiv versions upon my suggestions as a math reviewer, and we didn't have any discussion about that.
See e.g. arXiv:1609.06573v3,   arXiv:1209.1170v4 and a forthcoming paper by D. Gugnin.
A less positive example is arXiv:1512.05164v6 +  arXiv:1808.08363v2.

Best Regards, Arkadiy.

\smallskip
{\it (November 5)} Dear Adam, Dear Tye,

Attached please find an update of my note ( = extended version of math review).
Outside Remark 5 changes are minimal.

On October, 14 and 20 I asked you to update your interesting paper arxiv:1803.01765
adding missing details of the proof (they are described in my letters and in my note whose project I sent you on October 30).
On October, 22 you refused to do so.
Please note that I already have delayed math review on your paper.
So, unless I find an arxiv update of arxiv:1803.01765 on November, 8, that day I'll publish my note, and on November, 10 I'll submit the review itself.
Please let me know if making missing details publicly available requires more time than I expected, so that you need more time to prepare the arxiv update.

If your reliability standards are different from mine and you disagree that the proof in arxiv:1803.01765 is incomplete
(in the sense of arxiv:1702.04259v1, text after Conjecture 2), I encourage you to publicly state this.
I am willing to publish your statement in an update of my note.
(Of course, I will publish it literally, without any corrections or omissions.)
Using a non-trivial result without explicit statement (and so without a proof or a reference) is a common flaw.
So it is important for math community to have a common standard regarding such flaws, by having
(potentially different) public opinions on particular cases.

Best Regards, Arkadiy.

\smallskip
{\it (November 6)} Dear Arkadiy,

As we have stated previously, we completely stand by the content of our paper, which was carefully reviewed by a referee with expertise in the field and found to be correct.
We do not intend to make any changes to our arXiv posting beyond the published version, as you have not pointed out any true gaps or errors.
We recognize that our paper assumes a certain amount of familiarity with Heegaard Floer homology on the part of the reader,
but we feel that the level of detail is entirely within the standards of the field.
We have explained each mathematical question that you have asked, and therefore we do not intend to engage with you further about this matter.

We also wanted to notify you that we have reported this entire exchange to the editors of Mathematical Reviews for an impartial judgment.

Sincerely, Adam Levine and Tye Lidman

\smallskip
{\it (November 17)} Dear Adam, Dear Tye,

Thank you for making our discussion public by reporting my letters to the editors of MR.
The more mathematicians and math journals publicly state their reliability standards by considering particular examples,
the more competent could be decisions of math community, sponsors and tax-payers to support this or other trend in mathematics.
(The same applies to refereeing standards, see below.)
Even without your efforts it is in the competence of MR to reject a report or to balance it with a report from another mathematician.

Thank you also for stating, however implicitly, that {\it reliability standards of yours (and of the field) involve considering a proof complete even if
some parts of the proof are missing, some results are used without proof or reference, or even without explicit statement, as in the above example.}
Here  I call a proof {\it incomplete} if one mathematician should be able to expect from another

(1) to wait for another (`complete') proof {\it before} using results having such a proof;

(2) to recommend, as a referee, a revision (based on specific comments) {\it before} recommending publication of  results having such a proof;

(3) to work more on such a proof (in particular, send the text privately to a small number of mathematicians working on related problems), {\it before} submitting the text to a refereed journal or to arxiv.

(Unfortunately, shorter formulations of this notion were found to be potentially misleading.)

You stated your reliability standards only implicitly
(in spite of my encouragement to make an explicit statement in the letter of November 5).
You use the unexplained notion `true gaps or errors' instead of the notion of `incomplete proof' explained in practical terms (the above
explanation is cited from reference in my letter of November, 5).
Mathematician often avoid explicit statements and explanations in practical terms, when they do understand that
such a clarity could damage their reputation.
So it would be nice if you could either publicly explicitly state the statement above in italics, or update your paper arxiv:1803.01765
adding missing details of the proof.
Correcting flaws commands respect, while unwillingness to add details and correct mistakes is more dangerous to reliability of mathematics than mistakes and missing details themselves.

Thank you also for showing, however implicitly, your refereeing standards, by calling `careful' the report to your paper
(which report presumably contained no suggestions of improvement which you did not fulfill in the published version).
Let me explicitly show my refereeing standards (which I consider common).
If I do not have time to referee carefully enough (to show that some parts of the proof are missing, some results are used without proof or reference, or even without explicit statement, and to indicate that making a paper well-structured would make it significantly shorter or clearer),
then I would either refuse to referee or would call my report `a quick opinion' and suggest to the editors another referee.

Mathematical papers are written for the math community, not for Arkadiy Skopenkov.
Publication of an argument involves more responsibility and more efforts than writing a letter to a colleague describing this argument.
So your explanations to me alone give nothing to the community.
However, I hope that these explanations would not be hard to elaborate and incorporate into update of your paper arxiv:1803.01765.

In the attached update of my note please find Remark 4b, which concerns other relevant question than reliability.

Sincerely, Arkadiy.

\bigskip

[LL] A.S. Levine and T. Lidman, Simply connected, spineless 4-manifolds, Forum of Math., Sigma, 7 (2019) e14, 1--11, arxiv:1803.01765.

\end{document}